\documentclass[12pt]{article}
\usepackage{amsmath}
\usepackage{amssymb}
\usepackage{graphicx}
\usepackage{caption2}
\usepackage{amsfonts}
\usepackage{cite}
\usepackage{lineno}
\usepackage{color}
\usepackage{multirow}
\usepackage{subfigure}
\usepackage{bm}

\newcommand{\be}{\begin{equation}}
\newcommand{\ee}{\end{equation}}
\newcommand{\bea}{\begin{eqnarray}}
\newcommand{\eea}{\end{eqnarray}}
\newcommand{\bef}{\begin{figure}}
\newcommand{\ef}{\end{figure}}
\newcommand{\bt}{\begin{tabular}}
\newcommand{\et}{\end{tabular}}
\newcommand{\bno}{\begin{enumerate}}
\newcommand{\eno}{\end{enumerate}}

\def\3{\ss}

\catcode`\"=\active
\def"{\accent'177}

\begin{document}

\begin{center}
{\bf\large Non-uniqueness of smooth solutions of the Navier-Stokes equations from almost the same initial conditions}

Shijun Liao\footnote{Corresponding author, Email: sjliao@sjtu.edu.cn} and Shijie Qin\\

State Key Laboratory of Ocean Engineering, Shanghai, China\\
School of Ocean and Civil Engineering, Shanghai Jiaotong University, China
\end{center}

\begin{abstract}
 Using clean numerical simulation (CNS) which can give very accurate spatiotemporal trajectory of Navier-Stokes turbulence in a finite but long enough interval of time, we  give some numerical evidences that the  Navier-Stokes equations admit distinct global solutions from almost the same initial conditions whose difference is very small, i.e. even at  the order $10^{-40}$ of magnitude.    Hopefully  these  examples could provide some enlightenments for the uniqueness and existence of Navier-Stokes equations, which are related to one Millennium Prize Problem of Clay Institute.  
\end{abstract}

\hspace{-0.4cm}{\bf Keyword} Navier-Stokes equations, non-uniqueness, clean numerical simulation

\section{Introduction}

It has been widely accepted that turbulent flows can be described by the Navier-Stokes (NS) equations.  The turbulence is so important that the uniqueness and existence of smooth solution of NS equations is listed as one of Millennium Prize Problems~\cite{MillenniumProblem}.  Currently, Coiculescu and Palasek~\cite{Coiculescu2025-IM} proved that there exists the same initial data in the critical space $BMO^{-1}$ from which there exist two distinct global solutions of the NS equations, both smooth for all  $t>0$.  This is an elegant work, although their construction uses a simplified version of the mechanism \cite{Palasek2025-IM}  that  requires the data to have infinite energy.  

Can we explicitly give some numerical examples about non-uniqueness of solution of NS equations, which are very close to exact solution but with finite energy?   

Unfortunately,  this is impossible by means of traditional numerical methods.  Lorenz, who rediscovered the sensitivity dependence  on initial condition of chaos~\cite{lorenz1963deterministic}, discovered that numerical errors could have significant influences on global properties of chaotic systems~\cite{Lorenz2006Tellus}, because numerical noises as a kind of small artificial disturbance exponentially increase to macro-level due to the so-called ``butterfly-effect''.   Many researchers~\cite{Deissler1986PoF,  Vassilicos2023JFM,  boffetta2001predictability, berera2018chaotic, boffetta2017chaos} pointed out that turbulence is chaotic.  However, for all numerical algorithms including the direct numerical simulation (DNS) for  NS equations, numerical noises as a kind of small artificial disturbance  are unavoidable.  Thus, logically speaking, numerical noises of DNS of NS equations  should quickly increase to macro-level, which certainly  lead to distinct departure from exact trajectory of NS equations.    In other words, DNS of NS turbulence (i.e. turbulence  governed by NS equations) should be quickly polluted by numerical noises badly.  

In 2009,  Liao~\cite{Liao2009} proposed the so-called ``clean numerical simulation'' (CNS) so as to overcome the restrictions of the traditional algorithms.   Unlike DNS, numerical noises of CNS can be rigorously negligible in a finite   time-interval, which is long enough for statistics, by means of reducing both truncation error and  round-off error to a required small enough level~\cite{Liao2023book,Hu2020JCP, Qin2022JFM, Qin2024JOES, Liao-2025-JFM-NEC}.  In other words, CNS can give very accurate spatiotemporal trajectory of NS equations in a finite but long enough interval of time, which are very close to the exact solution.   Using CNS, the so-called ``noise-expansion cascade'' of turbulence was discovered~\cite{Liao-2025-JFM-NEC}, which indicates that all small distances of NS turbulence at different orders of magnitudes would quickly increase, one by one in order, to macro-level, and each of them could completely change the global properties of turbulence.             

In this paper, using CNS, we explicitly give two types of numerical examples of the NS equations, whose initial conditions are almost the same, but their corresponding  solutions are distinctly different not only in spatiotemporal trajectory and  spatial symmetry  but also even in statistics.   Hopefully, these explicit examples might give some enlightenments (or even better,  a roadmap) for mathematicians to prove non-uniqueness of smooth solution (with finite energy) of NS equations for $t>0$.           

\section{Numerical settings}

\subsection{Governing equation}

Consider a two-dimensional (2D) incompressible Kolmogorov flow \cite{obukhov1983kolmogorov, chandler2013invariant, wu2021quadratic} in a square domain  with a periodic boundary condition  under the so-called Kolmogorov forcing, which is governed by the  dimensionless Navier-Stokes equation in the form of stream function 
\begin{equation}
 \frac{\partial}{\partial t}\Big(\nabla^{2}\psi\Big)+\frac{\partial(\psi,\nabla^{2}\psi)}{\partial(x,y)}-\frac{1}{Re}\nabla^{4}\psi+n_K\cos(n_Ky)=0,       \label{eq_psi}
\end{equation}
where $\psi$ is the stream function,  $x,y\in[0,2\pi]$ are horizontal and vertical coordinates, $t$ denotes the time, $Re$ is the Reynolds number,   $n_K$ describes the Kolmogorov forcing scale,  
\[
 \frac{\partial(a,b)}{\partial(x,y)}=\frac{\partial a}{\partial x}\frac{\partial b}{\partial y}-\frac{\partial b}{\partial x}\frac{\partial a}{\partial y}      
\]
is the Jacobi operator,  $\nabla^{2}$ is the Laplace operator,  and  $\nabla^{4}=\nabla^{2}\nabla^{2}$,  respectively. The stream function $\psi$ satisfies the periodic boundary condition
\begin{equation}
\psi(x, y, t)=\psi(x\pm2\pi, y, t)=\psi(x, y\pm2\pi, t).       \label{boundary_condition}
\end{equation}
In order to investigate a state of turbulent flow,  we choose $n_K=16$ and $Re=2000$ thereafter in this paper.

\subsection{Initial conditions}

Let us consider the following two initial conditions
\begin{eqnarray}
 \psi_1(x,y,0) & = & -\frac{1}{2}\big[\cos(x+y) + \cos(x-y)],       \label{initial_condition1} \\
 \psi_2(x,y,0) & = & \psi_1(x,y,0) + \delta \; \sin(x+y),       \label{initial_condition2} 
\end{eqnarray}
where $\delta$ is a nonzero constant.   Let us choice here 
$\delta=10^{-10}$, $10^{-20}$, and $10^{-40}$,   i.e. $\delta\to 0$, corresponding to  the initial conditions
\begin{eqnarray}
 \psi_2(x,y,0) & = &   \psi_1(x,y,0) + 10^{-10}\sin(x+y),       \label{initial_condition2-1}\\ 
 \psi_2(x,y,0)  &= & \psi_1(x,y,0) + 10^{-20}\sin(x+y),       \label{initial_condition2-2}  \\
 \psi_2(x,y,0)  &=& \psi_1(x,y,0) + 10^{-40}\sin(x+y),       \label{initial_condition2-3} 
\end{eqnarray}
respectively.   Note that 
\begin{equation}
\left| \psi_2(x,y,0)  -  \psi_1(x,y,0) \right| = \left| \delta  \; \sin(x+y) \right| \leq \left| \delta \right|. 
\end{equation}
Thus, the three initial conditions (\ref{initial_condition2-1})-(\ref{initial_condition2-3}) are {\em almost} the same as $\psi_1(x,y,0)$.  For the sake of simplicity, we name the CNS results subject to the above-mentioned initial conditions (\ref{initial_condition1}), (\ref{initial_condition2-1}), (\ref{initial_condition2-2}), and (\ref{initial_condition2-3}), respectively,  as follows: Flow~CNS, Flow~CNS$'_1$, Flow~CNS$'_2$, and Flow~CNS$'_3$.  

It should be emphasized that  the initial condition $\psi_1(x,y,0)$ defined by (\ref{initial_condition1}) contains the  spatial symmetry of rotation and translation
\begin{equation}
\left\{
\begin{array}{ll}
\mbox{rotation}&: \;\; \psi(x,y,t)=\psi(2\pi-x,2\pi-y,t),\\
 \mbox{translation}&: \;\; \psi(x,y,t) =\psi(x+\pi,y+\pi,t).
 \end{array}  
 \right.     \label{symmetry_psi:A}
\end{equation}
However,   the initial condition $\psi_2(x,y,0)$ defined by (\ref{initial_condition2}) for arbitrary $\delta$ has only the spatial symmetry of translation 
\begin{equation}
 \psi(x,y,t)=\psi(x+\pi,y+\pi,t),        \label{symmetry_psi:B}
\end{equation}
although, for very small $\delta$,  it looks like having the same spatial symmetry as (\ref{symmetry_psi:A}).   

\subsection{Settings of CNS}

The above-mentioned NS equations are solved numerically by means of CNS in the square domain $[0,2\pi] \times [0,2\pi]$ on an uniform mesh $1024 \times 1024$ adopting the Fourier pseudo-spectral method for spatial approximation and the $M$th-order Taylor expansion with a time step $\Delta t = 10^{-3}$ for temporal approximation. Especially, different from DNS, we use {\em multiple-precision} with $N_s$ significant digits for all physical/numerical variables and parameters so as to decrease the round-off error to a small enough level.   Furthermore, to determine the so-called ``critical predictable time'' $T_c$,  another CNS result should be given by the same CNS algorithm but with even smaller numerical noise, i.e. using even larger $M$ and $N_{s}$ than those used previously, which confirms (by comparison) that the numerical noise of the former CNS result is  indeed rigorously negligible throughout the whole interval $t\in[0,T_c]$ so that it can be regarded as a  ``convergent''  benchmark solution that is very closed to its exact solution of the NS turbulence.  In addition, the self-adaptive CNS strategy \cite{Qin2023AAMM} and parallel computing are adopted to  increase the computational efficiency of the CNS algorithm  \cite{Liao2023book, Liao-2025-JFM-NEC}.   We  found that, to gain a reliable CNS in $t\in[0,300]$, $M=140$ \& $N_{s}=260$ are needed for Flow~CNS,  and $M=70$ \& $N_{s}=130$ are needed for Flow~CNS$'_1$, Flow~CNS$'_1$ and Flow~CNS$'_3$, respectively.  The corresponding CNS code can be downloaded via Github ({\color{blue}{https://github.com/sjtu-liao/2D-Kolmogorov-turbulence}}).

It should be emphasized that  all of our CNS results are carefully checked:  in the whole time interval $t\in[0,300]$, the grid spacing is always less than the averaged Kolmogorov scale and enstrophy dissipative scale~\cite{pope2001turbulent},  and the time-step is always small  enough, say, satisfying the Courant-Friedrichs-Lewy condition (i.e. Courant number $<1$)~\cite{courant1928partiellen}.

\section{Results}

\begin{figure}
    \begin{center}
        \begin{tabular}{cc}
             \includegraphics[width=2.5in]{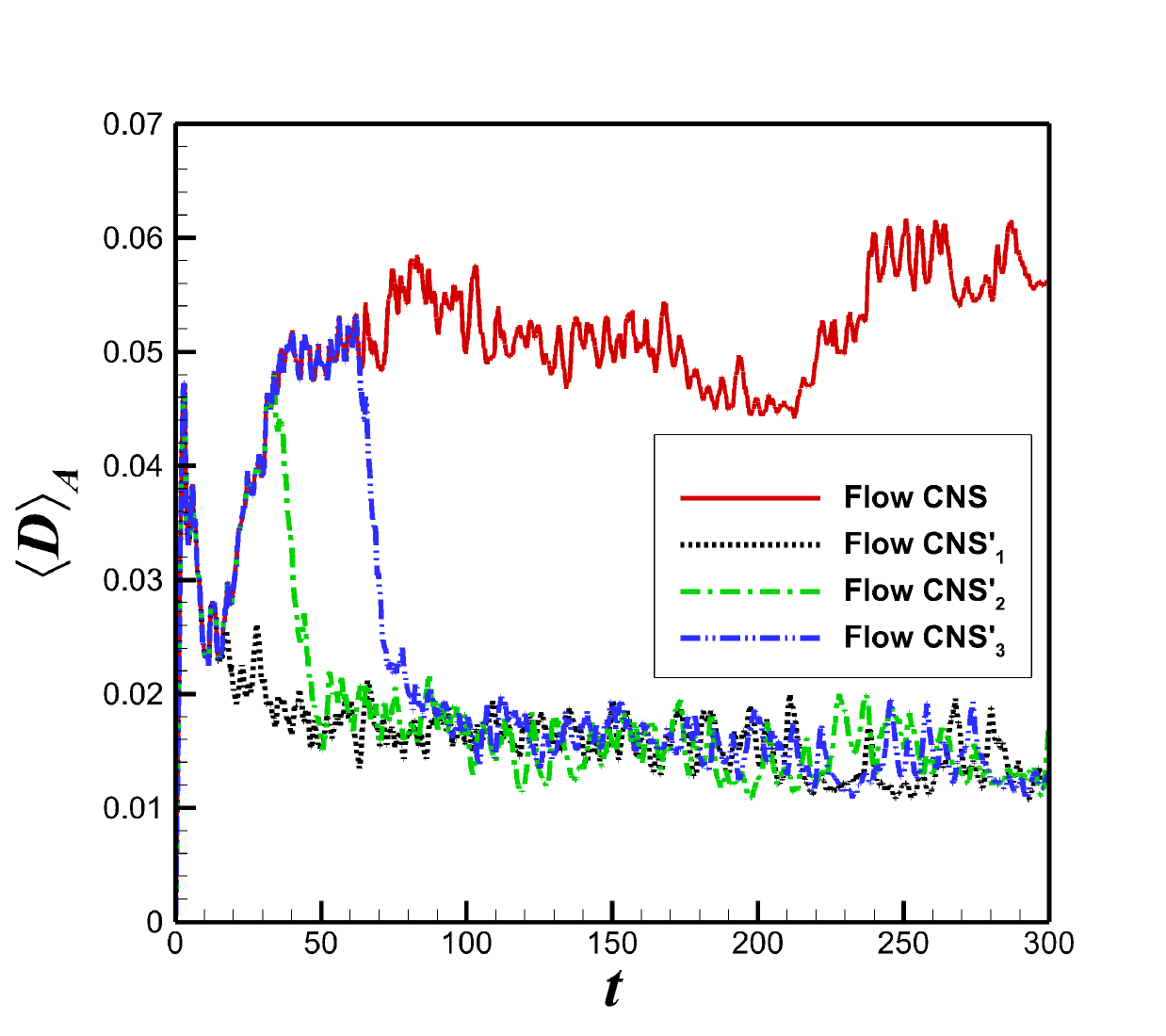}
        \end{tabular}
    \caption{Comparison of time histories of the spatially averaged kinetic energy dissipation rate $\langle D \rangle_A$ of the 2D turbulent Kolmogorov flow, governed by Eqs.~(\ref{eq_psi}) and (\ref{boundary_condition}) for $n_K=16$ and $Re=2000$, given by Flow~CNS (red solid line), Flow~CNS$'_1$ (black dash line), Flow~CNS$'_2$ (green dash-dot line), and Flow~CNS$'_3$ (blue dash-dot-dot line).}     \label{Fig:D_t}
    \end{center}
\end{figure}

\begin{figure}
    \begin{center}
        \begin{tabular}{cc}
             \includegraphics[width=2.0in]{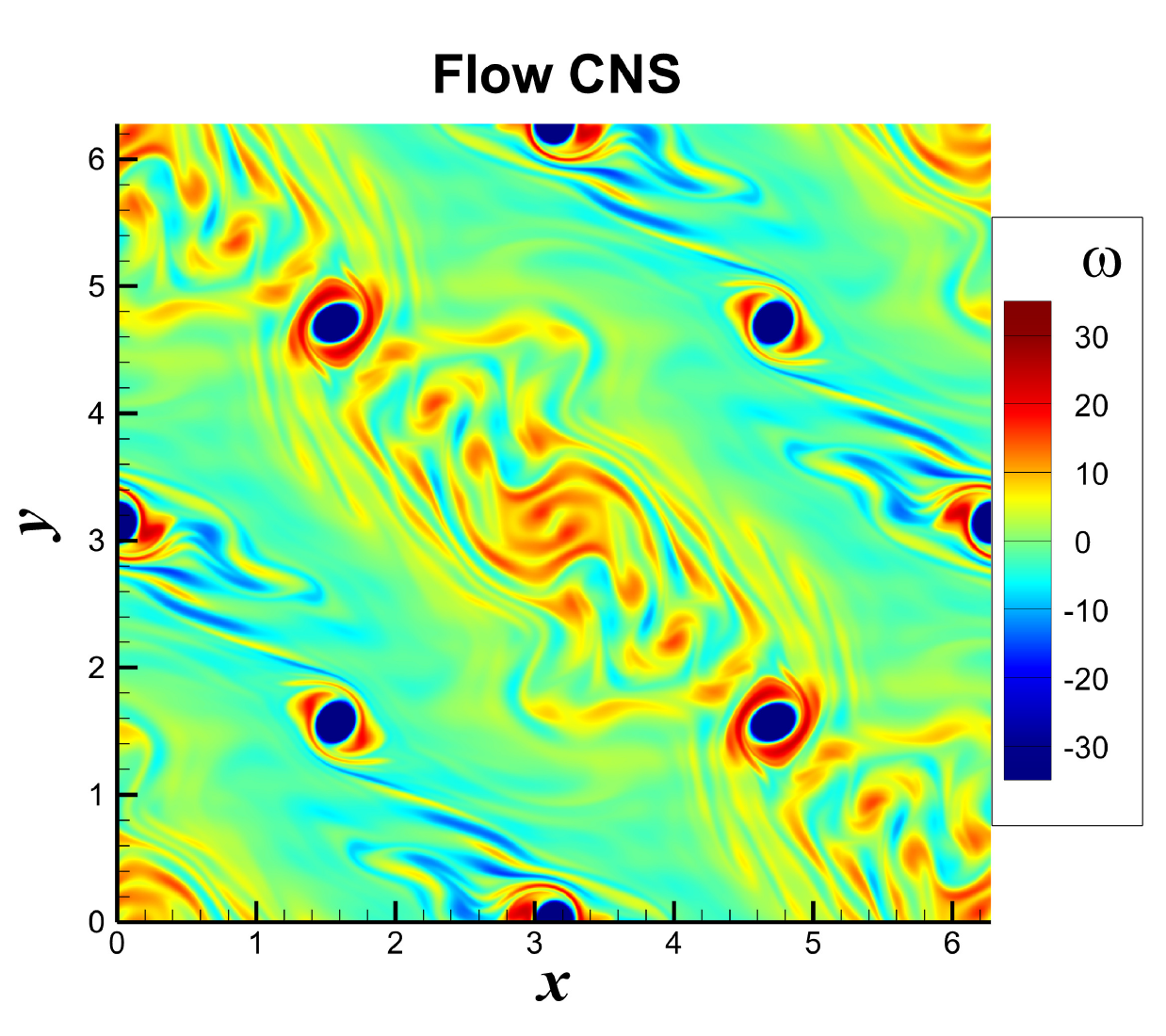}
             \includegraphics[width=2.0in]{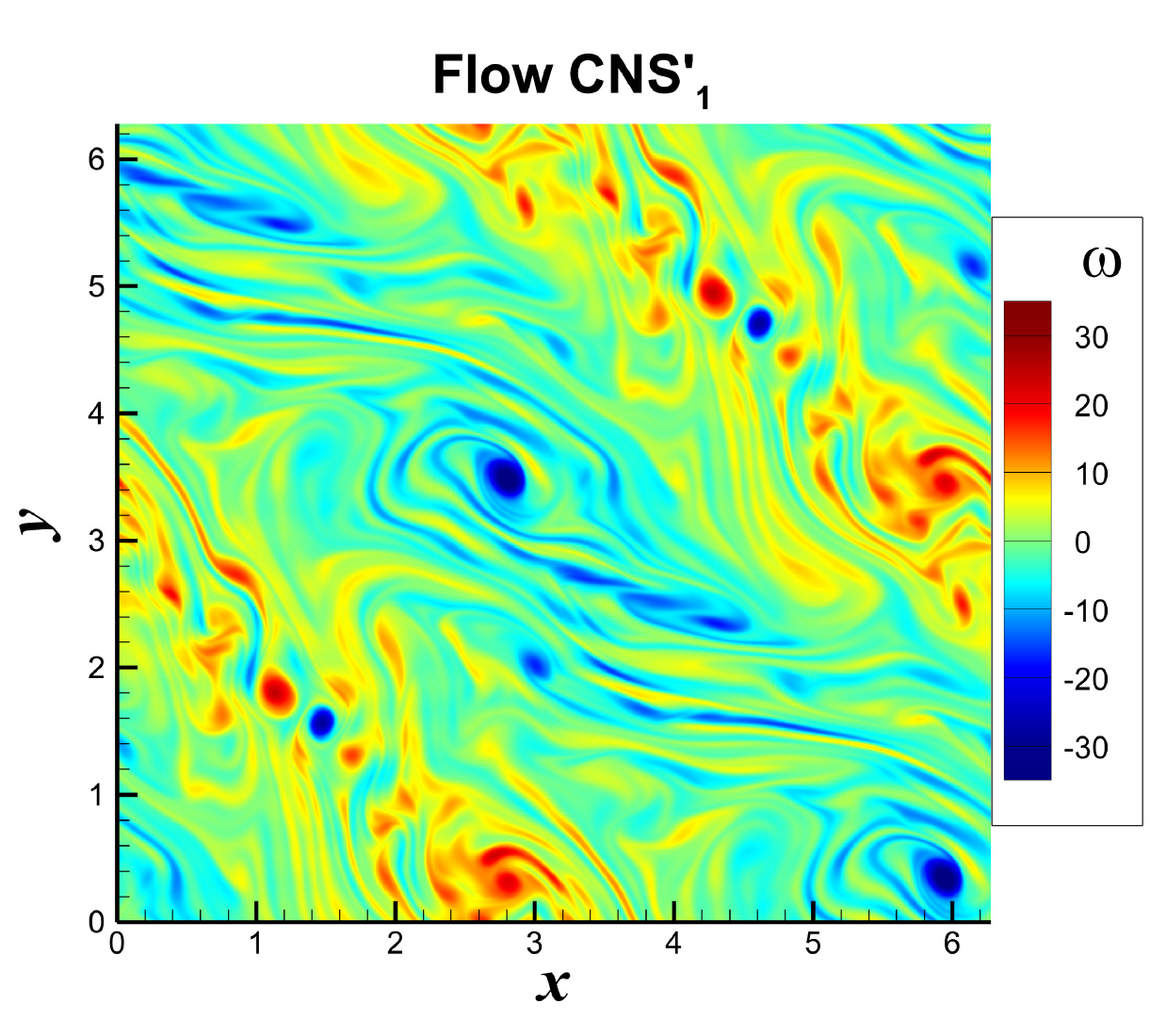}\\
             \includegraphics[width=2.0in]{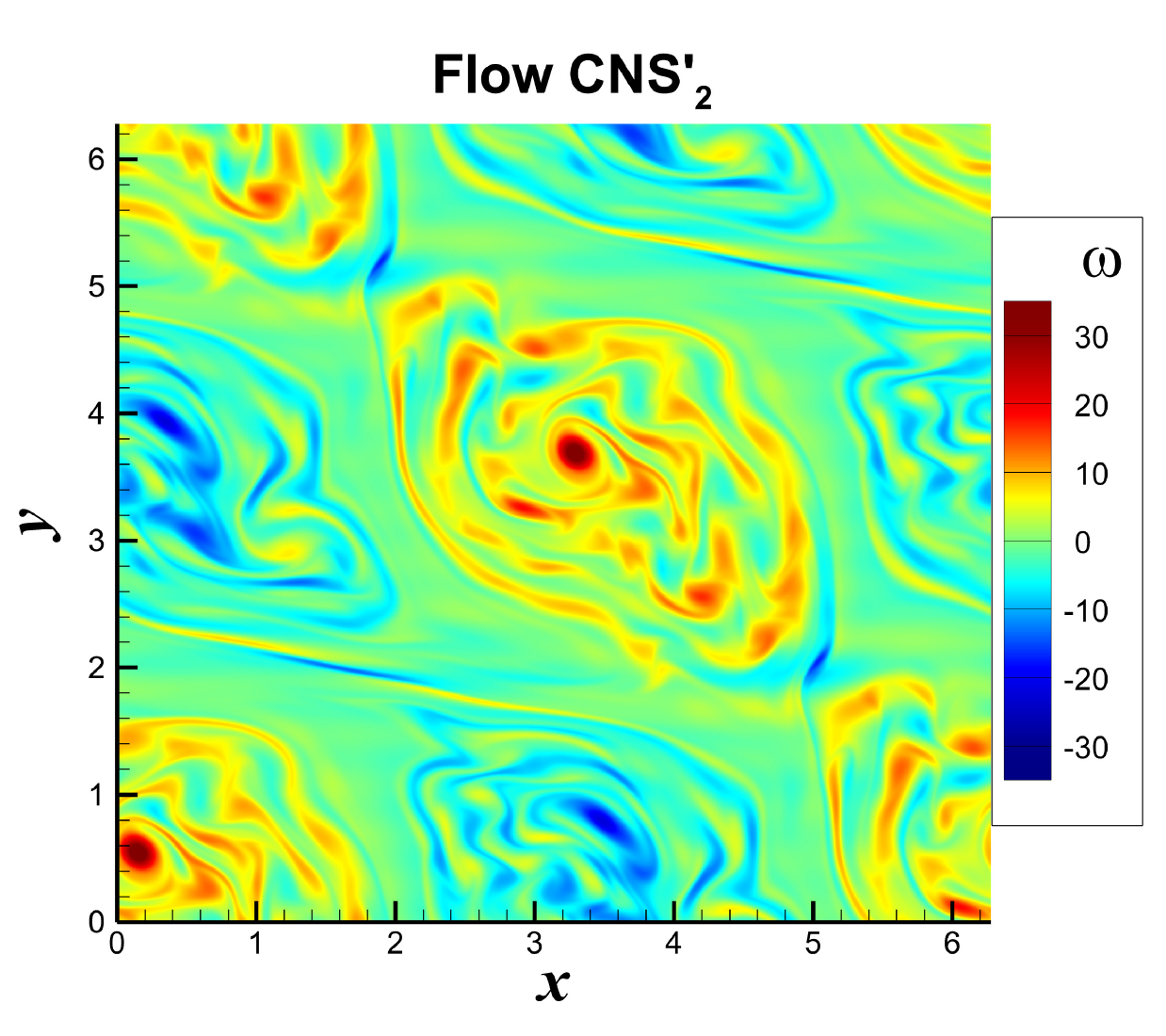}
             \includegraphics[width=2.0in]{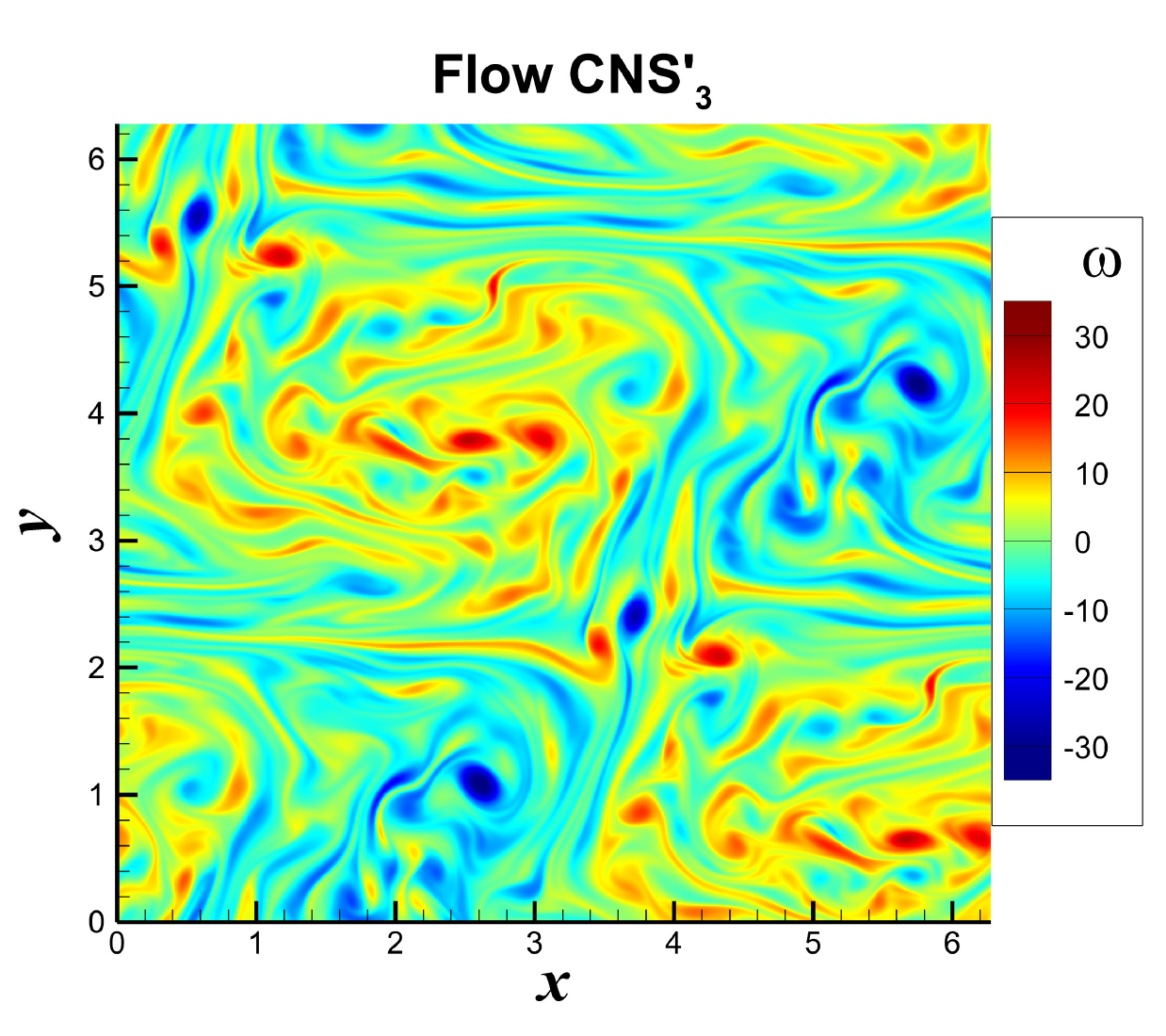}
        \end{tabular}
    \caption{Vorticity fields $\omega(x,y)$ at $t=200$ of the 2D turbulent Kolmogorov flow governed by (\ref{eq_psi}) and (\ref{boundary_condition}) for $n_K=16$ and $Re=2000$, given by CNS subject to the initial conditions (\ref{initial_condition1}) (marked by Flow~CNS), (\ref{initial_condition2-1}) (marked by Flow~CNS$'_1$), (\ref{initial_condition2-2}) (marked by Flow~CNS$'_2$), and (\ref{initial_condition2-3}) (marked by Flow~CNS$'_3$).
}     \label{Fig:Contour}
    \end{center}
\end{figure}

\begin{figure}
    \begin{center}
        \begin{tabular}{cc}
              \subfigure[]{\includegraphics[width=2.3in]{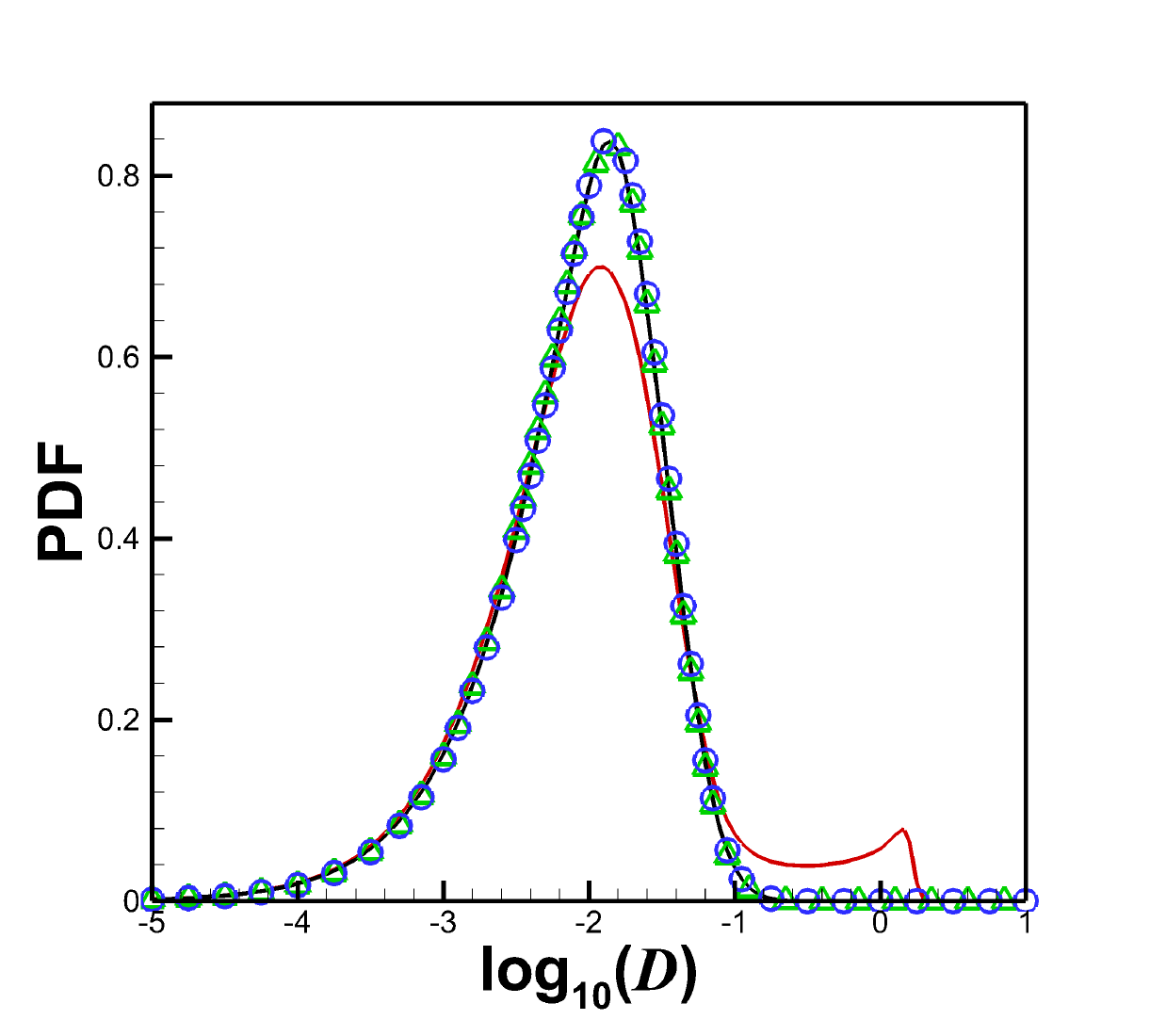}}
             \subfigure[]{\includegraphics[width=2.3in]{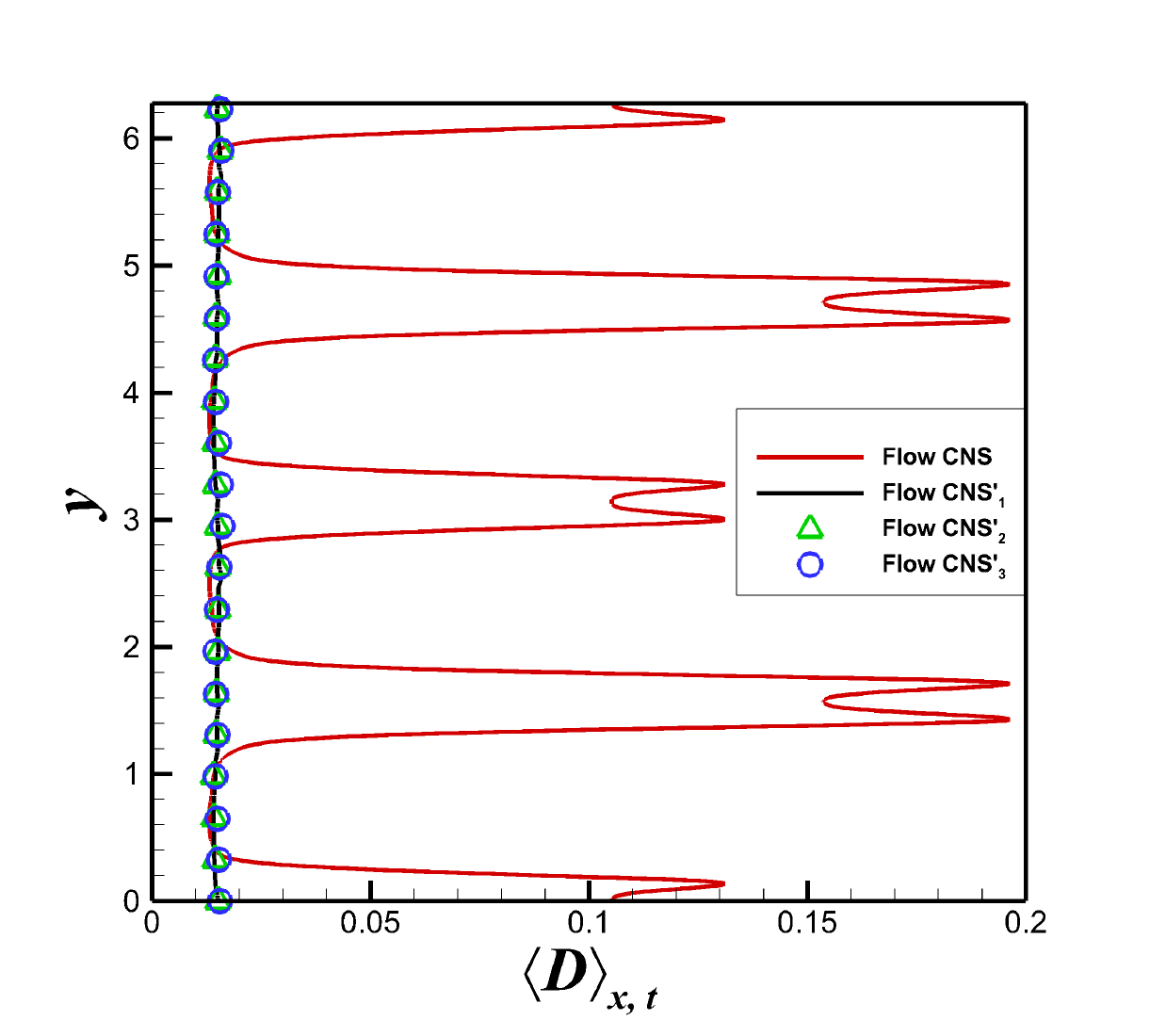}}
        \end{tabular}
    \caption{Comparison of (a) probability density functions (PDFs) of kinetic energy dissipation rate $D(x,y,t)$  and (b) spatio-temporal averaged  kinetic energy dissipation rate $\langle D\rangle_{x,t}(y)$ of the 2D turbulent Kolmogorov flow, governed by Eqs.~(\ref{eq_psi}) and (\ref{boundary_condition}) for $n_K=16$ and $Re=2000$, given by Flow~CNS (red line), Flow~CNS$'_1$ (black line), Flow~CNS$'_2$ (green line or triangle), and Flow~CNS$'_3$ (blue line or cirlce).}     \label{Fig:ED_y}
    \end{center}
\end{figure}

Note that the initial condition (\ref{initial_condition1}) is almost the same as another three initial conditions (\ref{initial_condition2-1})-(\ref{initial_condition2-3}).   It is found that at the beginning the Flow CNS$'_1$, Flow CNS$'_2$ and Flow CNS$'_3$ are almost the same as the Flow CNS, but after some durations  they depart from the Flow CNS one by one,  as shown in Fig.~\ref{Fig:D_t} for the spatially averaged kinetic energy dissipation rate $\langle D \rangle_A$.   Note that, when $t > 100$,  the value of $\langle D \rangle_A$ of the flow CNS is nearly two times larger than that of the Flow CNS$'_1$, Flow CNS$'_2$ and Flow CNS$'_3$.   
Besides,  the Flow CNS has the {\em same} spatial symmetry (\ref{symmetry_psi:A}) of rotation and translation as its initial condition  (\ref{initial_condition1}), but the Flow CNS$'_1$, Flow CNS$'_2$ and Flow CNS$'_3$ have the  {\em same} spatial symmetry (\ref{symmetry_psi:B}) of translation as its initial condition  (\ref{initial_condition2}), respectively, as shown in Fig.~\ref{Fig:Contour} at $t=200$, although these spatiotemporal trajectories are quite different.    This is mainly because, due to the chaotic characteristics of the  NS turbulence~\cite{Deissler1986PoF,  Vassilicos2023JFM,  boffetta2001predictability, berera2018chaotic, boffetta2017chaos},  the small disturbance $\delta\; \sin(x+y)$ in (\ref{initial_condition2})  increases quickly to a macro-level that destroys the spatial symmetry  (\ref{symmetry_psi:A}).  
In addition, the statistics of the Flow CNS$'_1$, Flow CNS$'_2$ and Flow CNS$'_3$ are the same but  
 have distinct  differences from those of the Flow CNS, as shown in Fig.~\ref{Fig:ED_y}.   Note that all of these CNS solutions of the NS turbulence have a finite energy density.      

All of these provide us a few examples that the NS equations could have distinct global solutions from {\em almost} the same initial conditions (\ref{initial_condition1}) and (\ref{initial_condition2}) for arbitrarily small $\delta$ or as $\delta \to 0$.   

\section{Concluding remarks and discussions}

Note that all of our numerical solutions of the NS turbulence are given by the CNS, whose numerical noises are much smaller than the corresponding exact solutions and thus are negligible in a finite  time-interval $t\in[0,300]$ that is long enough for statistics.  In other words, they are very close to the exact solution of the NS turbulence, and thus can give us some helpful enlightenments.   

First of all,  the different spatial symmetries are used to express the distinct difference of solutions of the NS turbulence.  Secondly,  the ``noise-expansion cascade''  characteristics of NS turbulence as a chaotic system is used so that a small disturbance such as $\delta \sin(x+y)$ in (\ref{initial_condition2})  can increase quickly to a macro-level so as to destroy the original spatial symmetry.   Although our CNS  results are numerical in a finite time interval $t\in[0,300]$ and rather special, they can be regarded as counter-examples for uniqueness of NS equations.   

Let $U_1({\bf r})$ denote an initial condition with a spatial symmetry,   $U_2({\bf r})=U_1({\bf r})+ \delta \; U_0({\bf r})$ is another initial condition with a different spatial symmetry (or without spatial symmetry), where ${\bf r}$ is a spatial vector and $\delta > 0$ is a real number, respectively.  Our CNS results mentioned above highly suggest   the following conjecture:  
   
{\em  
The Navier-Stokes equations admit distinct smooth global solutions for all $t>0$  from almost the same initial conditions  $U_1({\bf r})$  and  $U_2({\bf r})$ with $\left| U_1({\bf r})-U_2({\bf r})\right|\leq \delta$  for arbitrarily small $\delta$ (or as $\delta \to 0$). 
}    


\section*{Acknowledgements}
The calculations were performed on ``Tianhe New Generation Supercomputer'', National Supercomputer Center in Tianjin, China. This work is partly supported by National Natural Science Foundation of China (No. 12521002). 
	
\section*{Data and materials availability}
The corresponding CNS code can be downloaded via Github ({\color{blue}{https://github.com/sjtu-liao/2D-Kolmogorov-turbulence}}).  All data are available by sending requirement to the corresponding author.

\section*{Competing interests}
The authors declare no competing interests.

\bibliographystyle{elsarticle-num}
\bibliography{Kolmogorov}

\end{document}